\documentclass[12pt]{article}%
\usepackage{eurosym}
\usepackage{amsmath,amssymb,amsthm,amsfonts}
\usepackage{wasysym}
\usepackage{subcaption}
\usepackage{graphicx}
\usepackage{xcolor}
\usepackage[margin=1in]{geometry}
\usepackage{appendix}
\usepackage{textcomp}
\usepackage{listings}

\theoremstyle{definition}

\begin{document}

\title{Generalized Attracting Horseshoe in the R\"{o}ssler Attractor}
\author{Karthik Murthy\thanks{Bridgewater-Raritan High School}, Parth
Sojitra\thanks{Department of Electrical and Computer Engineering, New Jersey
Institute of Technology}, Aminur Rahman\thanks{Corresponding Author,
\url{ar276@njit.edu}} \thanks{Department of Mathematical Sciences, New Jersey
Institute of Technology} \thanks{Current Address: Department of Mathematics
and Statistics, Texas Tech University}, Ian Jordan\footnotemark[2]
\thanks{Current Address: Department of Applied Mathematics and Statistics,
Stony Brook University}, Denis Blackmore\footnotemark[4]}
\date{}
\maketitle

\begin{abstract}
We show that there is a mildly nonlinear three-dimensional system of ordinary
differential equations - realizable by a rather simple electronic circuit -
capable of producing a generalized attracting horseshoe map. A system
specifically designed to have a Poincar\'{e} section yielding the desired map
is described, but not pursued due to its complexity, which makes the
construction of a circuit realization exceedingly difficult. Instead, the
generalized attracting horseshoe and its trapping region is obtained by using
a carefully chosen Poincar\'{e} map of the R\"{o}ssler attractor. Novel
numerical techniques are employed to iterate the map of the trapping region to
approximate the chaotic strange attractor contained in the generalized
attracting horseshoe, and an electronic circuit is constructed to produce the
map. Several potential applications of the idea of a generalized attracting
horseshoe and a physical electronic circuit realization are proposed.

\end{abstract}


\section{Introduction}

\label{Sec: Intro}

The seminal work of Smale \cite{Smale1963} showed that the existence of a
\emph{horseshoe} structure in the iterate space of a diffeomorphism is enough
to prove it is chaotic. Often these diffeomorphisms arise from certain
Poincar\'{e} maps of continuous-time \emph{chaotic strange attractors} (CSA),
which in turn are discrete-time CSAs. Some examples of such attractors are the
Lorenz strange attractor \cite{Lorenz1963}, the R\"{o}ssler attractor
\cite{Rossler1976}, and the double scroll attractor \cite{Matsumoto84}. An
example of a Poincar\'{e} map of the Lorenz equations is the H\'{e}non map
\cite{Henon}, which can be further simplified to the Lozi map \cite{Lozi1978}.

In more recent years Joshi and Blackmore \cite{JB2014} developed an
\emph{attracting horseshoe} (AH) model for CSAs, which has two saddles and a
sink. This, however, negates the possibility of the H\'{e}non and Lozi maps,
which have two saddles. Fortunately the attracting horseshoe can be modified
into a \emph{generalized attracting horseshoe} (GAH), which can have either
one or two saddles while still being an attracting horseshoe \cite{JBR2017}.
This results in a quadrilateral trapping region. While extensive analysis was
done in Joshi \textit{et al. } \cite{JBR2017}, a simple concrete example
seemed to be illusive.

In this investigation we implement MATLAB codes to find the necessary
Poincar\'{e} map of the R\"{o}ssler attractor that would admit a quadrilateral
trapping region. The remainder of the paper is organized as follows; in Sec.
\ref{Sec: Algo} we give an overview of the algorithm with the MATLAB codes
relegated to the Appendix. Once we have the tools for our numerical
experiments, we first propose a carefully constructed GAH model in Sec.
\ref{Sec: Constructed}. Then, we give numerical examples of Poincar\'{e} maps
of the R\"{o}ssler attractor and the map of interest in Sec.
\ref{Sec: Poincare} and real world examples in Sec. \ref{Sec: Real World}.
Finally, we end discussions in Sec. \ref{Sec: Conclusion} with some concluding remarks.

\section{Poincar\'{e} map algorithm}

\label{Sec: Algo}

To produce a general Poincar\'{e} section of a flow we break up the program
into four parts: solving the ODE, computing a Poincar\'{e} section
perpendicular to either $x = 0$, $y = 0$, or $z = 0$, rotating the
Poincar\'{e} section, and iterating the Poincar\'{e} map. Solving the ODE is
standard through ODE45 on MATLAB, which executes a modified Runge-Kutta
scheme. Once we have our solution matrix we need to approximate the values of
first return maps from the discretized flow. By restricting the first return
onto a Poincar\'{e} section the iterate space of Poincar\'{e} map can be
visualized. This is easily done for a section perpendicular to the axes, but
in order to locate a highly specialized object such as the GAH we need to be
able to rotate the section. Once the desired section is found we can
experiment on iterating the points of trapping region candidates.

The first major task is approximating the first return map on a Poincar\'{e}
section of a flow. Much of the ideas of our initial first return map code came
from that of Didier Gonze \cite{PoincareCode}. Once the discretized flow is
found numerically a planar section for a certain value of $x$, $y$, or $z$ can
be defined, which in general will lie between pairs of simultaneous points.
Then we may draw a line between the pair through the planar section and
identify the intersecting point, which approximates a point of the first
return map. This can also be done with more simultaneous points in order to
get higher order approximations.

Once we can approximate a map for a section perpendicular to the axes we need
to have the ability to rotate and move the map to any position. This is where
our program completely diverges from that of \cite{PoincareCode}. While the
first instinct might be to try to rotate the section, it is equivalent to
rotate the flow in the opposite direction to the desired rotation of the
section. Once the flow is rotated, the code for the first return map can be
readily used. This gives us the ability to analyze the first return map of a
general Poincar\'{e} section.

Finally, we would like to not only compute a first return map, but also
compute the iterates of a Poincar\'{e} map of any system; that is, given an
initial condition on an arbitrary Poincar\'{e} section can we find the
subsequent iterates. To accomplish this, we solve the ODE for a given initial
condition on the planar section to find the first return. Once we have the
first return we record it's location and use that as the new initial
condition. This iterates the map for as many returns as desired, thereby
filling in a Poincar\'{e} map. Now we have the tools needed to run numerical
experiments on GAHs.

\section{A constructed GAH system}

\label{Sec: Constructed}

In this section, we give a brief description of the generalized attracting
horseshoe (GAH) map and devise a three-dimensional nonlinear ordinary
differential equation with a Poincar\'e section that produces it.

\subsection{The GAH map}

\begin{figure}[htbp]
\centering
\includegraphics[width=0.9\textwidth]{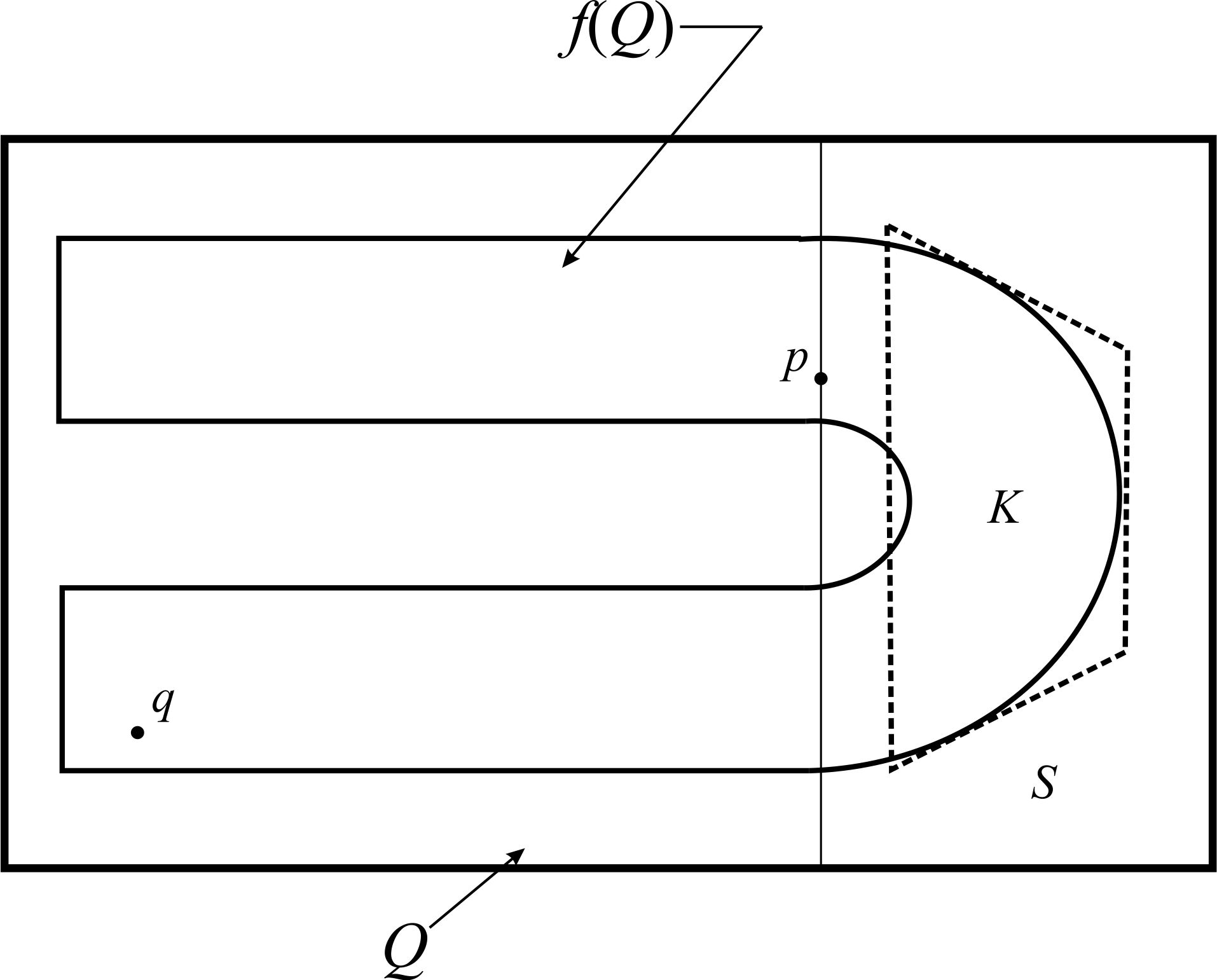}
\caption{A planar GAH with two saddle points}%
\label{fig:GAH}%
\end{figure}

The GAH is a modification of the AH that can be represented as a geometric
paradigm with either just one or two fixed points, both of which are saddles.
Figure \ref{fig:GAH} shows a rendering of a $C^{1}$ GAH with two saddle
points, which can be constructed as follows: The rectangle is first contracted
vertically by a factor $0<\lambda_{v}<1/2$, then expanded horizontally by a
factor $1<\lambda_{h}<2$ and then folded back into the usual horseshoe shape
in such a manner that the total height and width of the horseshoe do not
exceed the height and width, respectively of the trapping rectangle $Q$. Then
the horseshoe is translated horizontally so that it is completely contained in
$Q$. Obviously, the map $f$ defined by this construction is a smooth
diffeomorphism. Clearly, there are also many other ways to obtain this
geometrical configuration. For example, the map $f$ as described above is
orientation-preserving, and an orientation-reversing variant can be obtained
by composing it with a reflection in the horizontal axis of symmetry of the
rectangle, or by composing it with a reflection in the vertical axis of
symmetry followed by a composition with a half-turn. Another construction
method is to use the standard Smale horseshoe that starts with a rectangle,
followed by a horizontal composition with just the right scale factor or
factors to move the image of $Q$ into $Q$, while preserving the expansion and
contraction of the horseshoe along its length and width, respectively.

\begin{figure}[htbp]
\centering
\includegraphics[width=0.9\textwidth]{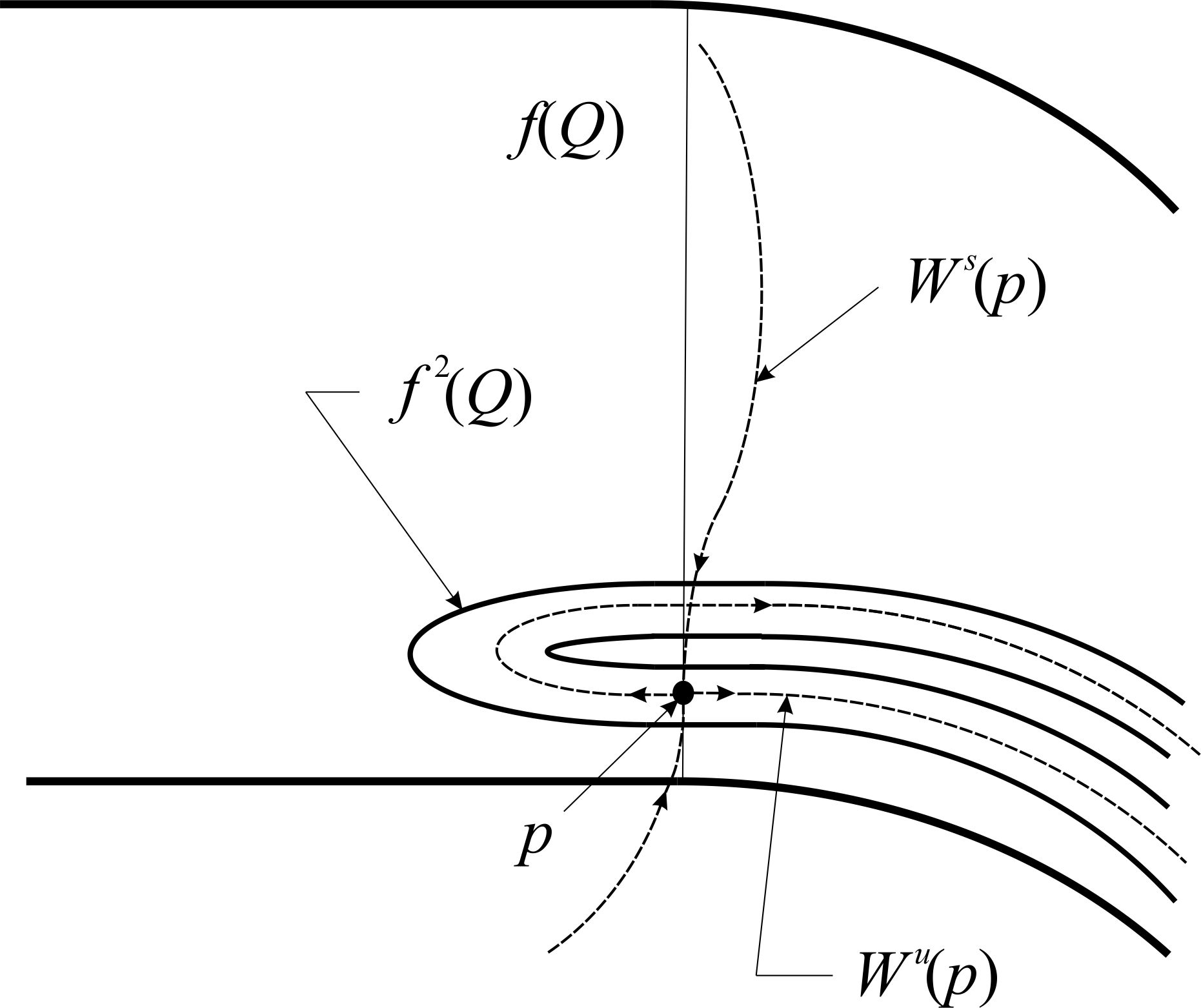}
\caption{Local (transverse) horseshoe structure of $f^{2}$ near $p$}%
\label{fig:f2p}%
\end{figure}

It is important to note that subrectangle $S$ with its left vertical edge
through $p$, which contains the arch of the horseshoe and the keystone region
$K$, plays a key role in the dynamics of the iterates of $f$. In particular,
we require that the map satisfy the following additional property, which is
illustrated in Fig. \ref{fig:f2p} :

\medskip

\begin{itemize}
\item[($\bigstar$)] \ $f$ \emph{maps the keystone region }$K$\emph{ }%
$($\emph{containing a portion of the arch of the horseshoe}$)$\emph{ to the
left of the fixed point }$p$ \emph{and the portion of its corresponding stable
manifold }$W^{s}(p)$ \emph{ containing }$p$\emph{ and contained in }$f(Q)$.
\end{itemize}

\medskip

\noindent The definition above and ($\bigstar$) can be shown to lead to the
conclusion that%

\[
\mathfrak{A}:=\overline{W^{u}(p)}={\displaystyle\bigcap\nolimits_{n=1}%
^{\infty}}f^{n}(Q),
\]
where $W^{u}(p)$ is the unstable manifold of $p$, is a global chaotic strange
attractor (CSA).

\medskip

The map above can be considered to be the paradigm for a GAH, but there are
many analogs. In fact, let $F:\tilde{Q}\rightarrow\tilde{Q}$ be any smooth
diffeomorphism of a quadrilateral trapping region $\tilde{Q}$ possessing a
horseshoe-like image with a keystone region $\tilde{K}$ containing a portion
of the arch of $F(\tilde{Q})$ analogous to that shown in Fig. \ref{fig:GAH}.
Suppose that the map is expanding by a scale factor uniformly greater than one
along the length of the horseshoe and contracting transverse to it by a scale
factor uniformly less than one-half in the complement of a subset of
$\tilde{Q}$ containing $\tilde{K}$. Then if $F$ satisfies an additional
property analogous to $(\bigstar)$, it maps $\tilde{K}$ into an open subset of
$\tilde{Q}$ to the left of the saddle point $\tilde{p}$, and%
\[
\mathfrak{A}:=\overline{W^{u}(\tilde{p})}={\displaystyle\bigcap\nolimits_{n=1}%
^{\infty}}F^{n}(\tilde{Q})
\]
\noindent is a global CSA.

\subsection{A GAH producing system}

We now construct an ODE in $\mathbb{R}^{3}$ with a Poincar\'{e} section that
is a GAH. The transversal we use is the following square in the $xz$-plane
defined in Cartesian and polar coordinates%
\begin{align}
Q_{0}  &  :=\{(x,y,z):0.05\leq x\leq1.05,y=0,-0.5\leq z\leq0.5\}\nonumber\\
\,  &  =\{(r,\theta,z):0.05\leq r\leq1.05,\theta=0,-0.5\leq z\leq0.5\}.
\label{e1}%
\end{align}
The trick is to find a relatively simple (necessarily nonlinear) $C^{1}$ ODE
having $Q_{0}$ as a transversal with an induced Poincar\'{e} first-return map
$P:Q_{0}\rightarrow Q_{0}\supset Q_{2\pi}:=P(Q_{0})$ such that $P(Q_{0}%
)\subset\mathrm{int}Q_{0}$ is a GAH. We chose the ODE based upon a rotation
about the $z$-axis so that the square evolves into the GAH as $Q_{0}$ makes a
full rotation. The first half of the metamorphosis takes care of the vertical
squeezing and horizontal stretching, while the second half produces the
folding. It is not difficult to show that the system (in cylindrical
coordinates)
\begin{equation}
\dot{r}=\frac{2\log(1.2)\sin^{2}\theta}{\pi}r,\;\dot{\theta}=1,\;\dot{z}%
=\frac{2\log5\sin^{2}\theta}{\pi}(z+0.2) \label{e2}%
\end{equation}
flows $Q_{0}$ to
\begin{equation}
Q_{\pi}:=\{(x,y,z):-1.26\leq x\leq-0.06,y=0,--0.26\leq z\leq-0.06\} \label{e3}%
\end{equation}
which is the original square in the radial half-plane plane corresponding to
$\theta=0$ stretched by a factor of 1.2 along the $x$-axis and squeezed by a
factor of 1/5 with respect to $z=-0.2$ along the $z$-axis in the radial half
plane corresponding to $\theta=\pi$. Consequently, (\ref{e2}) produces the
first half of the desired result comprising the stretching and squeezing for
$0\leq\theta\leq\pi$.

Note that (\ref{e2}) can be integrated directly to obtain the following for
$0\leq\theta\leq\pi$ and initial condition $\left(  r(0),\theta
(0),z(0)\right)  :=\left(  r_{0},\theta_{0},z_{0}\right)  :$%
\begin{align}
r(t)  &  =r(\theta)=r_{0}\left[  \frac{\log(1.2)}{2\pi}(2\theta-\sin
2\theta\right]  ,\nonumber\\
\theta(t)  &  =t,\label{e4}\\
z(t)  &  =z(\theta)=-0.2+(z_{0}+0.2)\left[  \frac{\log(0.2)}{2\pi}%
(2\theta-\sin2\theta\right]  .\nonumber
\end{align}
Now we have to attend to the folding for $\pi\leq\theta\leq2\pi$. For this we
use a rotation in planes orthogonal to a fixed circle in the $xy$-plane. In
these planes corresponding to a circle of radius $c$, given as $c=0.66$, we
define Euclidean coordinates with origin $r=0.66,z=0$ and corresponding polar
coordinates $(\rho,\phi)$ as%
\begin{equation}
\rho:=\sqrt{(r-0.66)^{2}+z^{2}}:=\sqrt{\tilde{r}^{2}+z^{2}}, \label{e5}%
\end{equation}
where $\tilde{r}:=r-0.66=\rho\cos\phi$ and $z:=\rho\sin\phi$. Then, when
$\pi\leq\theta\leq2\pi$, we take the folding part for $\phi\geq-\pi/2$ to be
\begin{equation}
\dot{\tilde{r}}=\dot{r}=-2\sin^{2}\theta\,\rho\sin\phi,\;\dot{\theta}%
=1,\;\dot{z}=2\sin^{2}\theta\,\rho\cos\phi, \label{e6}%
\end{equation}
or equivalently%
\begin{equation}
\dot{\tilde{r}}=\dot{r}=-2\sin^{2}\theta\,z,\;\dot{\theta}=1,\;\dot{z}%
=2\sin^{2}\theta\,\tilde{r}. \label{e7}%
\end{equation}
It is easy to verify from the above that $\rho$ is constant (call it $\rho
_{0}$) for the solutions of (\ref{e6}) or (\ref{e7}) and that the solution
initially (at $t=\theta=\pi$) satisfying $(\rho,\phi)=(\rho_{0},\phi_{0})$ is%
\begin{equation}
\tilde{r}=\tilde{r}(t)=\rho_{0}\cos\left(  \phi(t)+\phi_{0}\right)
,\;\theta=\theta(t)=t,\;z=z(t)=\rho_{0}\sin\left(  \phi(t)+\phi_{0}\right)  ,
\label{e8}%
\end{equation}
where%
\begin{equation}
\phi(t):=(t-\pi)-\sin t\cos t. \label{e9}%
\end{equation}

The above ((\ref{e6}) or (\ref{e7})) describes the folding field for $\pi
\leq\theta\leq2\pi$ and $-\pi/2\leq\phi$. In order to smoothly fill in the
rest of the field, we shall use the function%
\begin{equation}
\psi(\tilde{r}):=\left\{
\begin{array}
[c]{cc}%
0, & \tilde{r}\leq-0.6\\
\sin^{2}\left[  \frac{\pi}{1.2}(\tilde{r}+0.6)\right]  , & -0.6\leq\tilde
{r}\leq0
\end{array}
\right.  , \label{e10}%
\end{equation}
which can be recast as
\begin{equation}
\xi(r):=\left\{
\begin{array}
[c]{cc}%
0, & r\leq0.06\\
\sin^{2}\left[  \frac{\pi}{1.2}(r-0.06)\right]  , & 0.06\leq r\leq0.66
\end{array}
\right.  . \label{e11}%
\end{equation}

We have now assembled all the elements for defining an ODE that generates a
GAH Poincar\'{e}\ section. This ODE, which incorporates (\ref{e2}) and
(\ref{e7}) and is $\pi$-periodic in $\theta$, has the following form:
\begin{equation}
\dot{r}=R(r,\theta,z),\;\dot{\theta}=1,\;\dot{z}=Z(r,\theta,z), \label{e12}%
\end{equation}
subject to the initial condition%
\begin{equation}
\left(  r(0),\theta(0),z(0)\right)  =\left(  r_{0},0,z_{0}\right)  \in Q_{0},
\label{e13}%
\end{equation}
where%
\[
R:=\left\{
\begin{array}
[c]{cc}%
\frac{\log(1.2)\sigma(\theta)r}{\pi}, & 0\leq\theta\leq\pi\\
-\sigma(\theta)z, & (\pi\leq\theta\leq2\pi)\,\mathrm{and\,}\left(  \left(
(r\geq0.66)\,\mathrm{or\,}(z\geq0)\right)  =(-\pi/2\leq\phi\leq\pi)\right) \\
-\xi(r)\sigma(\theta)z, & (\pi\leq\theta\leq2\pi)\,\mathrm{and\,}\left(
(r<0.66)\,\mathrm{and\,}(z\in\lbrack-0.26,-0.06])\right) \\
0, & (\pi\leq\theta\leq2\pi)\,\mathrm{and\,}\left(  (r<0.66)\,\mathrm{and\,}%
(z<0)\,\mathrm{and\,}(z\notin\lbrack-0.26,-0.06])\right)
\end{array}
\right.  ,
\]%
\[
Z:=\left\{
\begin{array}
[c]{cc}%
\frac{(\log(.2)\sigma(\theta)(z+0.2)}{\pi}, & 0\leq\theta\leq\pi\\
\sigma(\theta)\tilde{r}, & (\pi\leq\theta\leq2\pi)\,\mathrm{and\,}\left(
\left(  (r\geq0.66)\,\mathrm{or\,}(z\geq0)\right)  =(-\pi/2\leq\phi\leq
\pi)\right) \\
0, & (\pi\leq\theta\leq2\pi)\,\mathrm{and\,}\left(  \left(
(r<0.66)\,\mathrm{and\,}(z<0)\right)  =(-\pi\leq\phi\leq-\pi/2)\right)
\end{array}
\right.  ,
\]
and%
\[
\sigma(\theta):=2\sin^{2}\theta=1-\cos2\theta.
\]
Finally, it is not difficult to show that the Poincar\'{e} section of the
transversal (and trapping region) $Q_{0}$ under the system (\ref{e12}) is a
GAH with image that is simply a 180-degree rotation of the horseshoe in
Fig. \ref{fig:GAH}. However, it appears that the construction of an electronic
circuit simulating (\ref{e12}) would be a rather formidable undertaking, so we
selected a simpler system; namely, the R\"{o}ssler attractor model, which is a
mildly nonlinear three-dimensional ODE that has a straightforward circuit realization.

\section{Poincar\'{e} maps and circuit realization of the R\"{o}ssler attractor}\label{Sec: Poincare}

\begin{figure}[htbp]
\centering
\includegraphics[width = 0.9\textwidth]{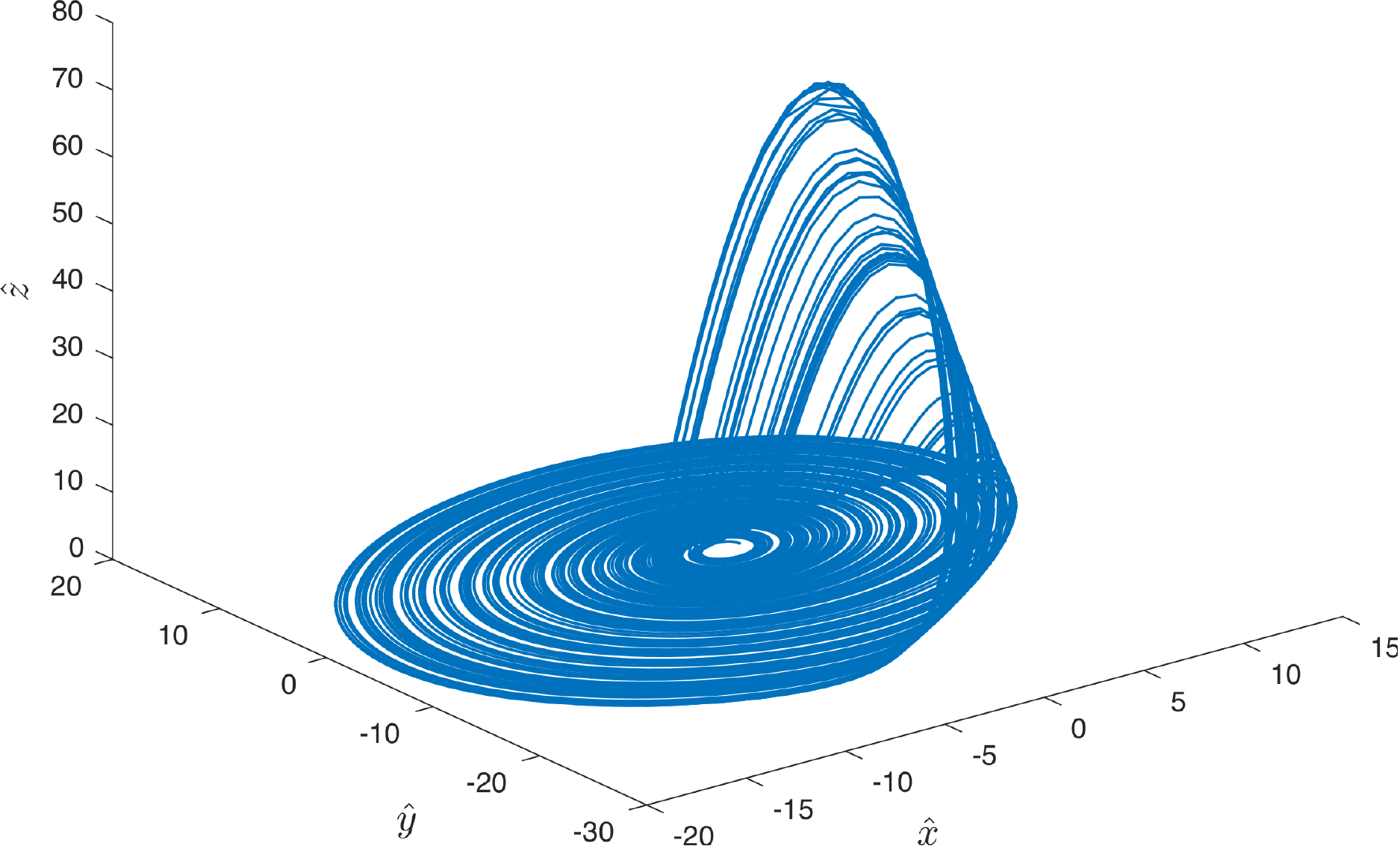}
\caption{The R\"{o}ssler attractor with parameters $a = 0.2$, $b = 0.1$, and $c = 10$, and a rotation
(represented by $\hat{x}$ and $\hat{y}$) of $\theta= 2\pi/5$ in spherical coordinates.}%
\label{Fig: RosslerAttractor}%
\end{figure}

We consider the R\"{o}ssler attractor
\begin{equation}%
\begin{split}
\dot{x}  &  = y - z\\
\dot{y}  &  = x + ay\\
\dot{z}  &  = b + z(x - c);
\end{split}
\label{Eq: Rossler}
\end{equation}
where we use the parameters $a = 0.2$, $b = 0.1$, and $c = 10$. This produces
the chaotic strange attractor in Fig. \ref{Fig: RosslerAttractor} and it can
also be realized by a rather simple electronic circuit.

\subsection{The Poincar\'{e} map}

\begin{figure}[htbp]
\centering
\includegraphics[width = 0.9\textwidth]{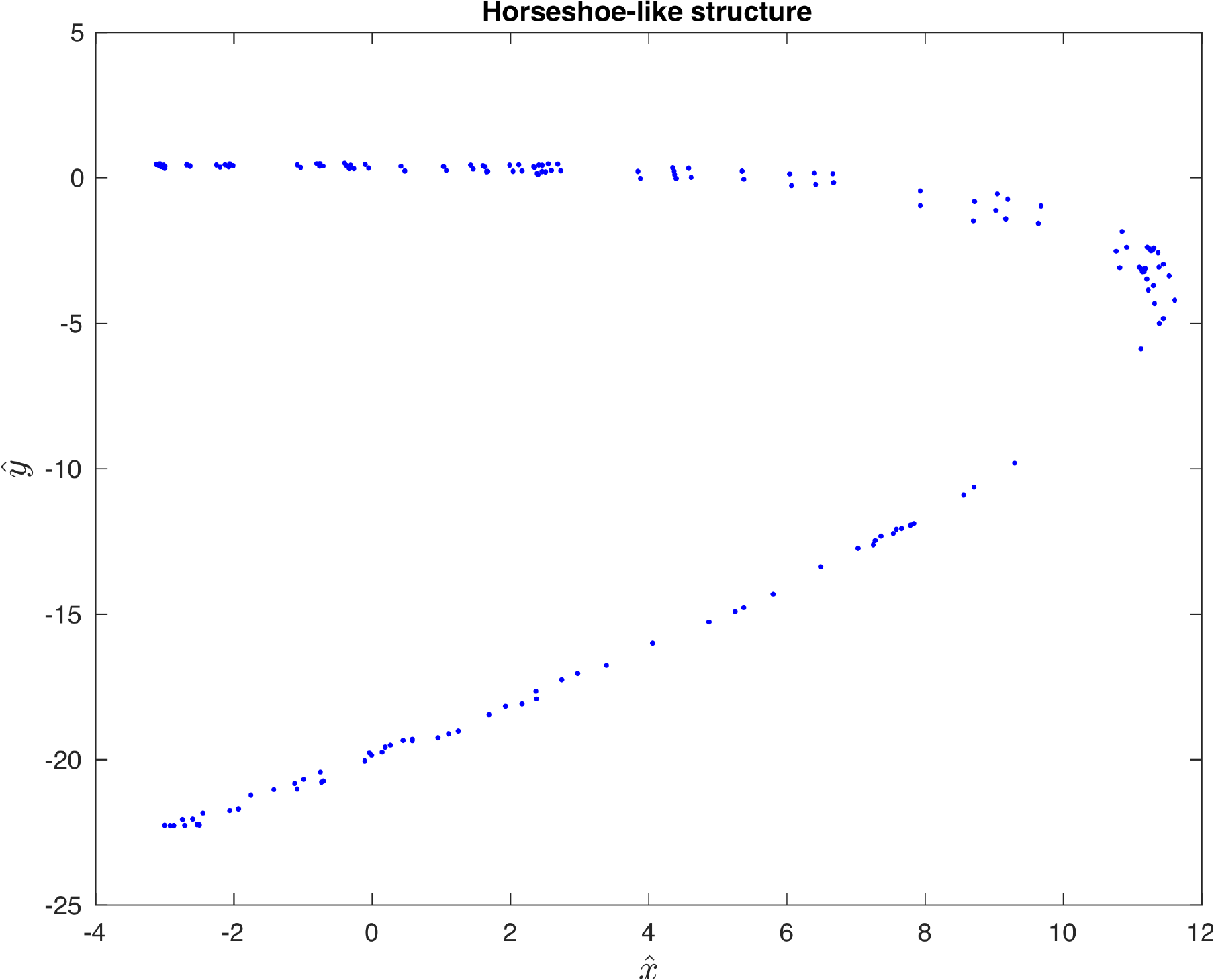}
\caption{Poincar\'{e} section ($r = 5, \theta= 2\pi/5$) of the R\"{o}ssler attractor containing a
horseshoe-like structure. Plot is shown in the rotated frame.}%
\label{Fig: Horseshoe}%
\end{figure}

One can use the algorithm in Sec. \ref{Sec: Algo} to compute any Poincar\'{e}
section of the attractor, however what we are particularly interested in is
identifying a trapping region for a generalized attracting horseshoe. Assuming
the system contains a GAH, we first look for a Poincar\'{e} section with a
horseshoe-like structure as shown in Fig. \ref{Fig: Horseshoe}.

\begin{figure}[htbp]
\centering
\includegraphics[width = 0.9\textwidth]{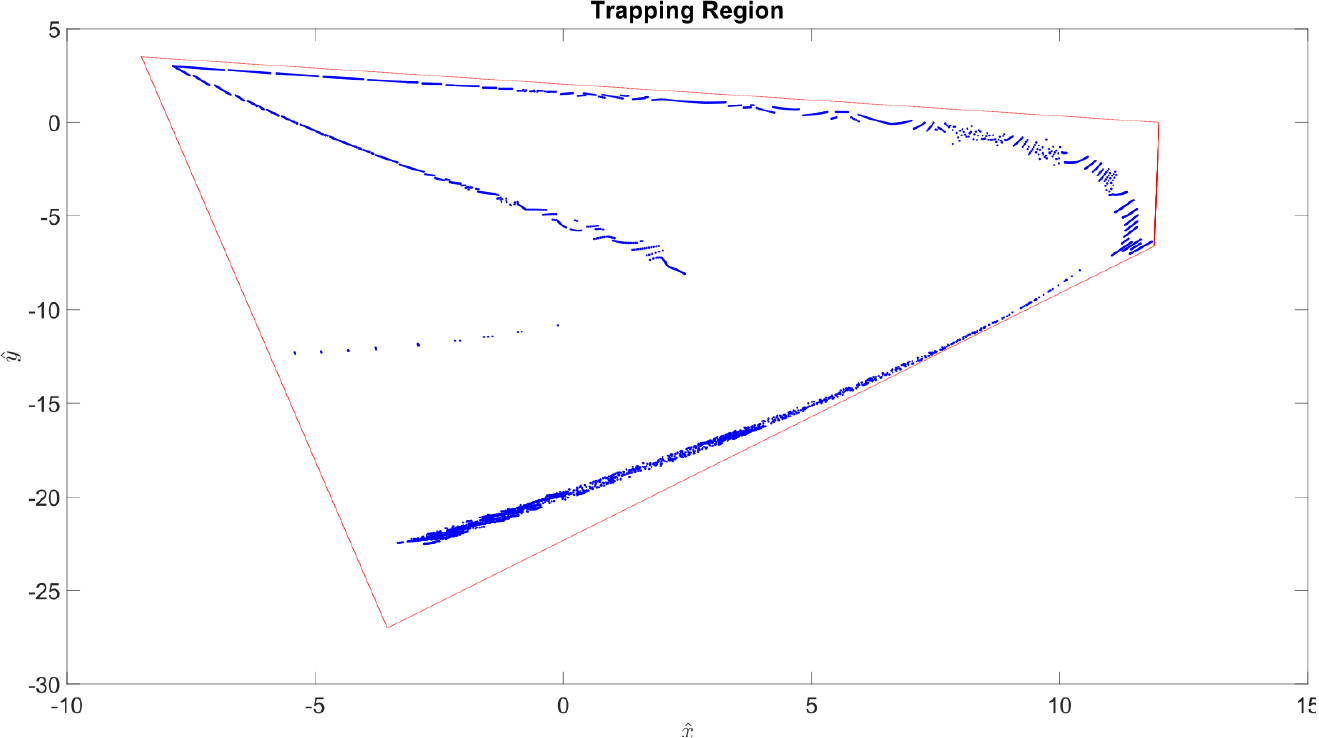}
\caption{The first return
(blue markers) of the quadrilateral trapping region (red markers) with
vertices located at $(\hat{x},\hat{y}) = (-3.55, -27),\, (11.91, -6.6),\, (12,
0),\, (-8.5, 3.5)$. While the quadrilateral edges look ``continuous'', it
should be noted that it is in fact discretized using four thousand points,
which are then mapped back to the Poincar\'{e} section ($r = 5, \theta=
2\pi/5$). Plot is shown in the rotated frame with $\hat{x}$ and $\hat{y}$
denoting rotated axes.}%
\label{Fig: TrappingRegion}%
\end{figure}

Now, if we can find a trapping region around this horseshoe, we will have
shown evidence for the existence of a GAH. First we identify vertices of a
quadrilateral that fully encompasses the horseshoe-like structure. Then using
a recursive algorithm (described in Sec. \ref{Sec: Algo}) we compute the first
return map of those vertices on that particular Poincar\'{e} section; i.e.,
the first iteration of the Poincar\'{e} map of those points. If the iterates
are contained within that quadrilateral, the points on the quadrilateral
itself can be tested. In Fig. \ref{Fig: TrappingRegion} four thousand points
on the quadrilateral are iterated and it is illustrated that this first return
is completely contained in the quadrilateral. While this is not a proof, the
grid spacing on the quadrilateral provides compelling evidence that this may
be a trapping region for the GAH. 

\begin{figure}[htbp]
\centering
\includegraphics[width = 0.9\textwidth]{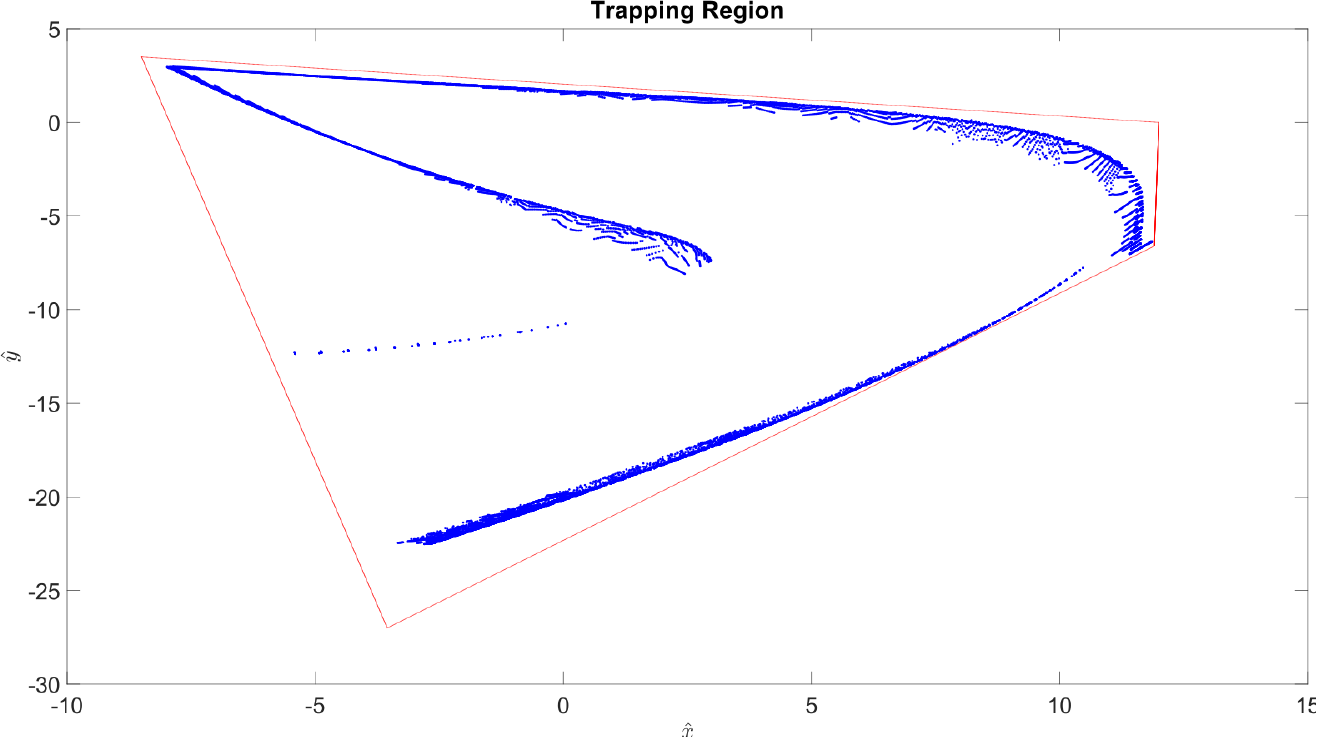}
\caption{First
five iterations of the Poincar\'{e} map (blue markers) of the quadrilateral
trapping region (red markers) with vertices located at $(\hat{x},\hat{y}) =
(-3.55, -27),\, (11.91, -6.6),\, (12, 0),\, (-8.5, 3.5)$. While the
quadrilateral edges look ``continuous'', it should be noted that it is in fact
discretized using four thousand points, which are then mapped back to the
Poincar\'{e} section ($r = 5, \theta= 2\pi/5$). Plot is shown in the rotated
frame with $\hat{x}$ and $\hat{y}$ denoting rotated axes.}%
\label{Fig: TrappingRegion5thOrder}%
\end{figure}

In order to provide more compelling evidence, we compute higher order
iterations of the Poincar\'{e} map in Fig. \ref{Fig: TrappingRegion5thOrder}.

\subsection{Circuit realization of the R\"{o}ssler system}

\begin{figure}[htbp]
\centering
\includegraphics[width=0.9\textwidth]{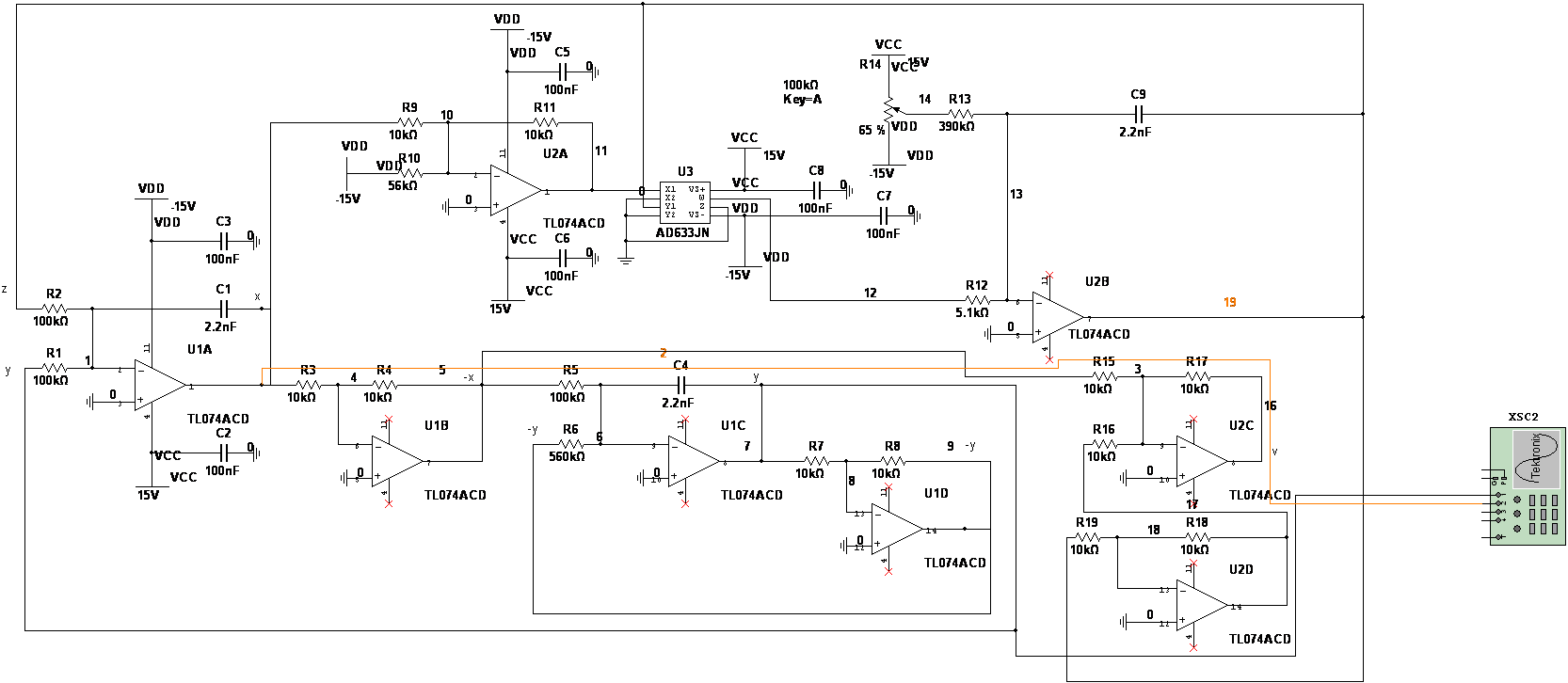}
\caption{Multisim circuit diagram for R\"{o}ossler attractor}
\label{fig:circuit}
\end{figure}

It happens that there are several known examples of electronic circuits
realizing the R\"{o}ssler attractor system. We chose the one, obtained from
\cite{RosslerCircuit}, shown in Fig. \ref{fig:circuit} with a list of components
in Table \ref{Tab: Parts}.

\begin{figure}[htbp]
\centering
\includegraphics[width = 0.32\textwidth]{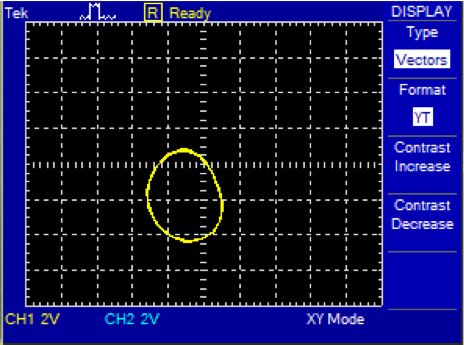}
\includegraphics[width = 0.32\textwidth]{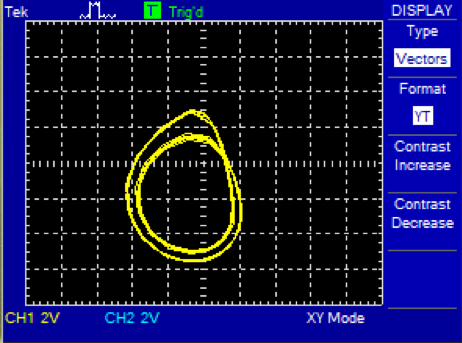}
\includegraphics[width = 0.32\textwidth]{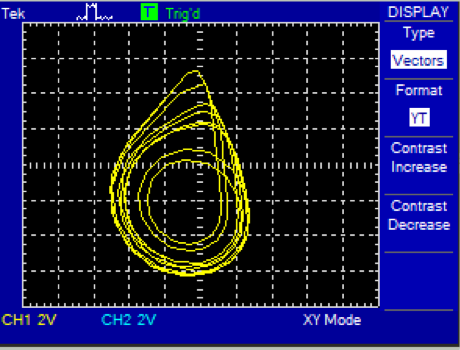}
\caption{Multisim outputs of the Rossler attractor showing a period doubling Hopf bifurcation leading to chaos.}
\label{Fig: Multisim}
\end{figure}

The physical realization of the Rossler attractor circuit was constructed using summing amplifiers, integrators, and a multipliers.  Due to the nature of this system, the operational amplifier must operate within $\pm 15$ volts in order to avoid clipping of the Rossler Attractor output waveform. In this circuit, resistors were used to represent constant values for parameters a and b in 
\eqref{Eq: Rossler}. A potentiometer was used to vary the parameter value of b in order to observe the bifurcations of the physical system.  We first test the circuit on \emph{Multisim} and observe the aforementioned bifurcations in Fig. \ref{Fig: Multisim}.

\begin{figure}[htbp]
\centering
\includegraphics[width=0.9\textwidth]{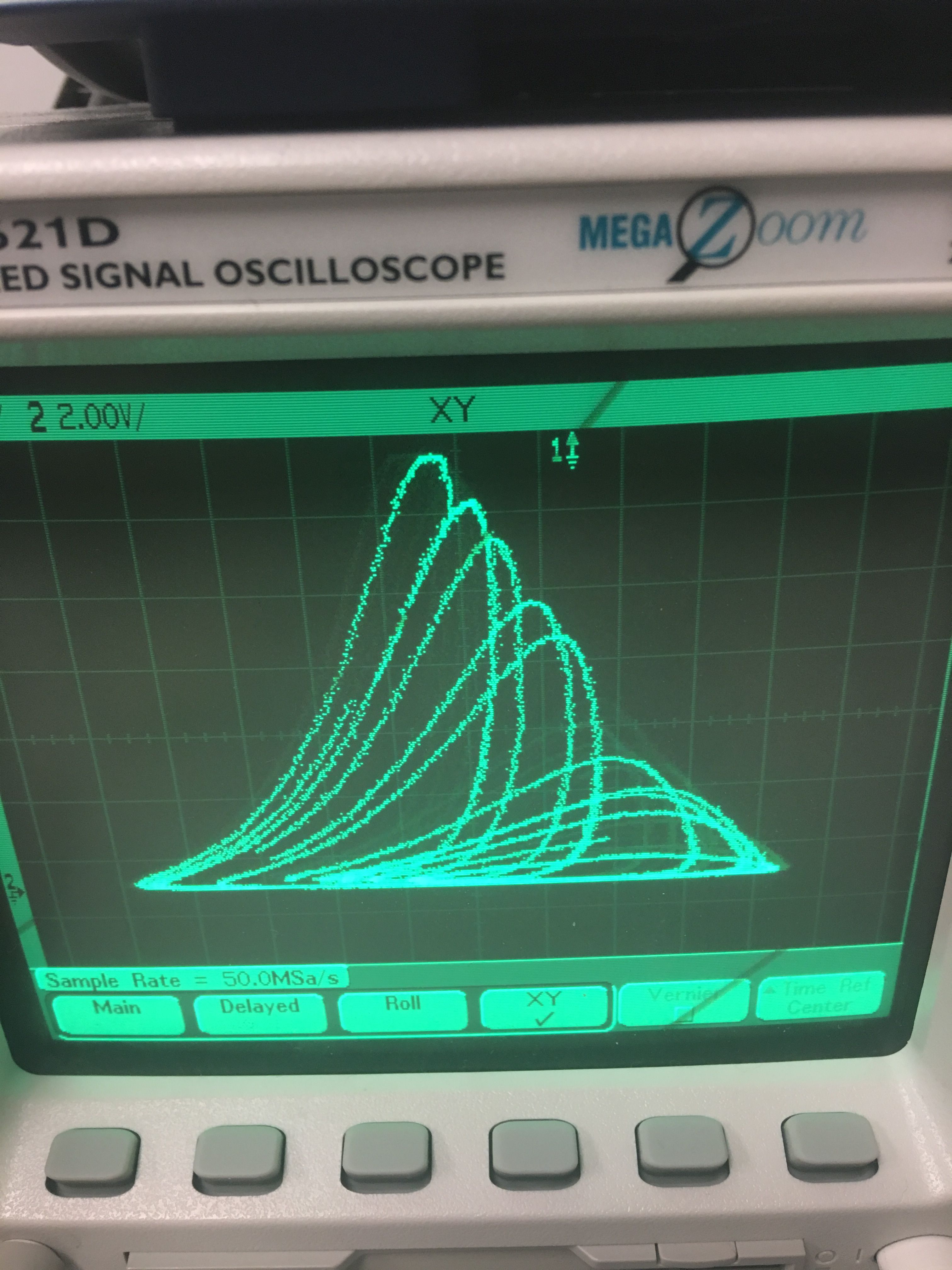}
\caption{Oscilloscope output from R\"{o}ssler attractor circuit}
\label{fig:oscil}%
\end{figure}

Next we built the circuit and observed oscilloscope outputs as shown in Fig. \ref{fig:oscil}. The
Poincar\'{e} section that we chose was a particular vertical plane through the
top arch of the output shown (see also Fig. \ref{Fig: RosslerAttractor}). The
acceptable planes were obtained by trial and error via varying the system
parameters and rotation of the plane about a vertical axis through the apex of
the arc.

\begin{table}[htbp]
\caption{List of components for the R\"{o}ossler attractor circuit}
\centering
\setlength{\tabcolsep}{12pt} 
\renewcommand{\arraystretch}{1}
\begin{tabular}{l||c|c|c|}
Type & Quantity & Code\\
\hline
$10 k\Omega$ Resistor & 11 &  \\
$100 k\Omega$ Resistor & 3 &  \\
$390 k\Omega$ Resistor & 1 &  \\
$56 k\Omega$ Resistor & 1 &  \\
$560 k\Omega$ Resistor & 1 &  \\
$5.1 k\Omega$ Resistor & 1 &  \\
$100 k\Omega$ Potentiometer & 1 &  \\
\hline
$100 nF$ Capacitor & 6 &  \\
$2.2 nF$ Capacitor & 3 &  \\
\hline
Op-Amp & 2 & AD633JN \\
Multiplier & 1 & TL074CN \\
\hline
\end{tabular}
\label{Tab: Parts}
\end{table}

\section{Potential applications}\label{Sec: Real World}

One can imagine several practical applications of devices containing
electronic circuit realizations of a GAH. Two, which are related to
communications and intelligence gathering, immediately come to mind: First,
the circuit could be embedded in a communication receiving device, and tuned
to certain "static" frequencies different from those in the expected incoming
messages. The strong global attracting characteristics of the circuit would
separate the static from the incoming messages, thereby enhancing the
receiving capabilities of system. In effect, the GAH circuit would filter out
the static.

Secondly, a stationary or compact mobile device incorporating the GAH circuit
could be used to penetrate and analyze various communication systems. Either
be connecting remotely in the case of a stationary device or directly for a
mobile version, the global attracting properties could be employed to extract
crucial characteristics of the system to which it is connected. Moreover, the
same attracting features of the GAH circuit device could be used to absorb
various parts of sent messages that would render them useless, false or somply misleading.

The two rather basic applications mentioned provide just a glimpse of the
possible applications of GAH circuits, most of which would probably be related
to information systems, data collection and filtering. Moreover, there are
more applications that could exploit the chaotic strange attractor associated
with a GAH circuit. For, example a GAH circuit device could be used either to
control chaos, introduce chaos or adjust the fractal dimension of outputs of a
variety of applicable processes based on dynamical systems.

\section{Conclusions}
\label{Sec: Conclusion}
We constructed a rather complicated nonlinear three-dimensional ordinary differential 
equation (ODE) having a Poincar\'e section that is a GAH map, but is not particularly 
amenable to electronic circuit realization, which was a goal of the investigation. 

So,instead of the initial ODE, we selected the R\"{o}ssler attractor; a mildly nonlinear 
three-dimensional ODE that has a reasonably simple circuit realization
and can actually produce GAH maps for carefully chosen Poincar\'e sections. We constructed 
the corresponding GAH circuit and used a novel iteration procedure to generate good 
approximations of the chaotic strange attractors associated to the GAH maps.

Finally, in addition to the experimental and analytic aspects of our investigation, we 
discussed a number of potential practical applications of the GAH circuit. Most of 
the envisioned applications were in the realms of communication and information gathering.

\section*{Acknowledgment}

The authors would like to thank the NJIT Provost Research Grants for funding
this research. K.M. was funded by the Provost high school internship and P.S.
was funded by Phase-1 Provost Undergraduate Research Grant, with D.B. as
faculty mentor and A.R. as graduate student mentor. K.M. appreciates the
support of Bridgewater-Raritan High School, P.S. and I.J. appreciate the
support of the Electrical and Computer Engineering Department at NJIT, and
A.R. and D.B. appreciate the support of the Department of Mathematical
Sciences and the Center for Applied Mathematics and Statistics (CAMS) at NJIT.

\bibliographystyle{unsrt}
\bibliography{GAH_ExpMath}

\begin{thebibliography}{10}

\bibitem{Smale1963}
Stephen Smale.
\newblock Diffeomorphisms with many periodic points.
\newblock In S.~Carins, editor, {\em Differential and Combinatorial Topology},
  pages 63--80, Princeton, NJ, 1963. Princeton University Press.

\bibitem{Lorenz1963}
Edward~Norton Lorenz.
\newblock Deterministic nonperiodic flow.
\newblock {\em J. Atoms. Sci.}, 20:130--141, 1963.

\bibitem{Rossler1976}
O.E. R\"{o}ssler.
\newblock An equation for continuous chaos.
\newblock {\em Phys. Lett. A}, 57A(5):397--398, 1976.

\bibitem{Matsumoto84}
T.~Matsumoto.
\newblock A chaotic attractor from {C}hua's circuit.
\newblock {\em IEEE Transactions on Circuits and Systems},
  CAS-31(12):1055--1058, December 1984.

\bibitem{Henon}
M.~H{\'e}non.
\newblock A two-dimensional mapping with a strange attractor.
\newblock {\em Commun. Math. Phys.}, 50(1):69--77, 1976.

\bibitem{Lozi1978}
R~Lozi.
\newblock Un attracteur \'{e}trange (?) du type attracteur de h\'{e}non.
\newblock {\em Le Journal de Physique Colloques}, 39(C5):9--10, 1978.

\bibitem{JB2014}
D.~Blackmore and Y.~Joshi.
\newblock Strange attractors for asymptotically zero maps.
\newblock {\em Chaos, Solitons and Fractals}, 68:123--138, 2014.

\bibitem{JBR2017}
Y.~Joshi, D.~Blackmore, and A.~Rahman.
\newblock Generalized attracting horseshoe and chaotic strange attractors.
\newblock (under review), 2017.

\bibitem{PoincareCode}
Didier Gonze.
\newblock http://homepages.ulb.ac.be/~dgonze/info/matlab/poincare.m, February
  2011.

\bibitem{RosslerCircuit}
Glen K.
\newblock http://www.glensstuff.com/rosslerattractor/rossler.htm.

\end{thebibliography}

\pagebreak

\begin{appendices}
\section{Poincar\'{e} map codes}
\subsection{Attractor data}
\begin{lstlisting}
clear
t=[0,1000];
xinit=[0 1 0];  % Random initial point
%Solving system using ODE45
[t,x]=ode45(@(t,x) ODESys(x), t, xinit);
%%% Rotating flow to pick out desired Poincare section %%%
q = pi/2.5;
xtrans = [1 0 0; 0 cos(q) -sin(q); 0 sin(q) cos(q)];
ytrans = [cos(q)  0 sin(q); 0 1 0; -sin(q) 0 cos(q)];
ztrans = [cos(q) -sin(q) 0 ; sin(q) cos(q) 0; 0 0 1];
newdata = x * ztrans;
newerdata = [t, newdata];
save 'myfile.dat' newerdata -ascii; % Saves flow in myfile.dat
plot3(newdata(:,1),newdata(:,2),newdata(:,3),'LineWidth',1);
alw = 0.75;    % AxesLineWidth
fsz = 14;      % Fontsize
xlabel('$$\hat{x}$$','Interpreter','Latex','fontsize',14)
ylabel('$$\hat{y}$$','Interpreter','Latex','fontsize',14)
zlabel('$$\hat{z}$$','Interpreter','Latex','fontsize',14)
\end{lstlisting}
\subsection{ODE system}
\begin{lstlisting}
function dx= ODESys(x)
%%%%%%%%%%%%%%%%%%%%%%%%%%%%%%%%%%%%%%%%%%%%%%%%%%%%%%
%%% Rossler Equation %%%
%%% Parameters %%%
a= .2;
b= .1;
c= 10;
%%% Rossler ODEs %%%
dx1=-x(2)-x(3);
dx2=x(1)+a*x(2);
dx3=b+x(3)*(x(1)-c);
dx=[dx1;dx2;dx3];
%%%%%%%%%%%%%%%%%%%%%%%%%%%%%%%%%%%%%%%%%%%%%%%%%%%%%%
\end{lstlisting}
\subsection{Horseshoe-like structure}
\begin{lstlisting}
D=load('myfile.dat');   % The flow is saved in myfile.dat
var=4;          % Chooses the rotated coordinate axis for Poincare section
cut=5;          % Chooses the Poincare section in rotated space; e.g.
%   cut = 5 => the Poincare section is located at
%   r = 5, and theta = 2pi/5 D)
%%% Picks out other two variables %%%
for i = 2:4
if i ~= var
varo1=i;
i = i+1;
if i ~= var
varo2 = i;
else
varo2 = i+1;
end
break;
end
end
x=D(:,var);     % Extract time series for variable var
xt=x(2:end-1);  % Value of discretized flow on var axis
xp=x(1:end-2);  % Previous time value
xm=x(3:end);    % Next time value
%%% Finds which pairs of points are closest to the Poincare section %%%
k=find(((xp<cut) & (xt>=cut)) | ((xp>cut) & (xm<=cut)));
%%% Finds points on the flow intersecting with the Poincare section %%%
R=zeros(size(D));
for i=1:length(k)
R(i,:)=D(k(i),:);
end
plot(R(:,varo1),R(:,varo2),'b.')
% hold on
% plot([A(:,1); A(1,1)],[A(:,2); A(1,2)],'r.','MarkerSize',1)
% hold off
xlabel('$$\hat{x}$$','Interpreter','Latex','fontsize',14)
ylabel('$$\hat{y}$$','Interpreter','Latex','fontsize',14)
% title(sprintf('Trapping Region'))
title('Horseshoe-like structure')
alw = 0.75;    % AxesLineWidth
fsz = 14;      % Fontsize
\end{lstlisting}
\subsection{Poincar\'{e} map}
\begin{lstlisting}
function firstRet = PoincareMap(xinit) % xinit calls a particular
%    ordered triplet
t=[0,1000];
axis = 4;  % Chooses the rotated coordinate axis for Poincare section
cutter = 5;  % Chooses the Poincare section in rotated space; e.g.
%   cut = 5 => the Poincare section is located at
%   r = 5, and theta = 2pi/5
%%% Solving system using ODE45 %%%
q = -pi/2.5;  % Counter rotation to pick out the proper initial conditions
% in "ODESys"
ztrans = [cos(q) -sin(q) 0 ; sin(q) cos(q) 0; 0 0 1];
xinit = xinit * ztrans;
[t,x]=ode45(@(t,x) ODESys(x), t, xinit);
%%% Rotating flow to pick out desired Poincare section %%%
q = pi/2.5; % theta = 2pi/5
ztrans = [cos(q) -sin(q) 0 ; sin(q) cos(q) 0; 0 0 1];
newdata = x * ztrans;
newerdata = [t, newdata];
x= newerdata(: , axis);     % Extract time series for variable var
xt=x(2:end-1);  % Value of discretized flow on var axis
xp=x(1:end-2);  % Previous time value
xm=x(3:end);    % Next time value
m = length(xt);
for i = 1:m
if((xp(i, 1)<cutter) && (xt(i, 1)>=cutter)) || ((xp(i, 1)>cutter)...
&& (xm(i, 1)<=cutter))
firstRet = newerdata(i, :);
break
end
end
end
\end{lstlisting}
\subsection{Trapping region}
\begin{lstlisting}
D=load('myfile.dat');   % The flow is saved in myfile.dat
%%%  Vertices of quadrilateral %%%
% A = [-3.6 -27 5; 11.91 -6.5 5; 12 0 5; -8 3 5];  % This one works
A = [-3.55 -27 5; 11.91 -6.6 5; 12 0 5; -8.5 3.5 5];  % A better one
%%%  Points between the vertices on the quadrilateral %%%
x1 = [A(1,1):(A(2,1)-A(1,1))/1000:A(2,1)]';
y1 = (A(1,2)-A(2,2))/(A(1,1)-A(2,1))*(x1 - A(1,1)) + A(1,2);
z1 = ones(size(x1))*5;
x2 = [A(2,1):(A(3,1)-A(2,1))/1000:A(3,1)]';
y2 = (A(2,2)-A(3,2))/(A(2,1)-A(3,1))*(x2 - A(2,1)) + A(2,2);
z2 = ones(size(x2))*5;
x3 = [A(3,1):(A(4,1)-A(3,1))/1000:A(4,1)]';
y3 = (A(3,2)-A(4,2))/(A(3,1)-A(4,1))*(x3 - A(3,1)) + A(3,2);
z3 = ones(size(x3))*5;
x4 = [A(4,1):(A(1,1)-A(4,1))/1000:A(1,1)]';
y4 = (A(4,2)-A(1,2))/(A(4,1)-A(1,1))*(x4 - A(4,1)) + A(4,2);
z4 = ones(size(x4))*5;
A = [x1 y1 z1; x2 y2 z2; x3 y3 z3; x4 y4 z4];
var=4;          % Chooses the rotated coordinate axis for Poincare section
cut=5;          % Chooses the Poincare section in rotated space; e.g.
%   cut = 5 => the Poincare section is located at
%   r = 5, and theta = 2pi/5
%%% Picks out other two variables %%%
for i = 2:4
if i ~= var
varo1=i;
i = i+1;
if i ~= var
varo2 = i;
else
varo2 = i+1;
end
break;
end
end
x=D(:,var);     % Extract time series for variable var
xt=x(2:end-1);  % Value of discretized flow on var axis
xp=x(1:end-2);  % Previous time value
xm=x(3:end);    % Next time value
%%% Finds which pairs of points are closest to the Poincare section %%%
k=find(((xp<cut) & (xt>=cut)) | ((xp>cut) & (xm<=cut)));
%%% Computes iterates of A %%%
Q = zeros(length(A),4);
parfor i = 1:length(A)
Q(i,:) = PoincareMap(A(i, :));
end
%%% Manually add in number of higher order %%%
%%%     returns in order to use parfor %%%
%%% Second Returns %%%
Q1 = zeros(size(Q));
Q_spare = Q(:,2:4);
parfor i = 1:length(Q)
Q1(i,:) = PoincareMap(Q_spare(i, :));
end
%%% Third Returns %%%
Q2 = zeros(size(Q1));
Q_spare = Q1(:,2:4);
parfor i = 1:length(Q1)
Q2(i,:) = PoincareMap(Q_spare(i, :));
end
%%% Fourth Returns %%%
Q3 = zeros(size(Q2));
Q_spare = Q2(:,2:4);
parfor i = 1:length(Q)
Q3(i,:) = PoincareMap(Q_spare(i, :));
end
plot(Q(:,varo1),Q(:,varo2),'b.',Q1(:,varo1),Q1(:,varo2),'b.',...
Q2(:,varo1),Q2(:,varo2),'b.',Q3(:,varo1),Q3(:,varo2),'b.');
hold on
plot([A(:,1); A(1,1)],[A(:,2); A(1,2)],'r.','MarkerSize',1)
hold off
xlabel('$$\hat{x}$$','Interpreter','Latex','fontsize',14)
ylabel('$$\hat{y}$$','Interpreter','Latex','fontsize',14)
title(sprintf('Trapping Region'))
alw = 0.75;    % AxesLineWidth
fsz = 14;      % Fontsize
\end{lstlisting}
\end{appendices}

\end{document}